\documentclass[11pt,a4paper]{article}
\usepackage{graphicx} % avoid epsfig or earlier such packages
\usepackage{url}      % for URLs
\usepackage{amsmath}  % many want amsmath extensions
\newtheorem{thm}{Theorem}[section]
\newtheorem{cor}[thm]{Corollary}
\newtheorem{lem}[thm]{Lemma}

\newtheorem{rem}[thm]{Remark}
\DeclareMathOperator{\deriv}{d}
\DeclareMathOperator{\expon}{e}
\title{\bf Stochastic inequalities for single-server loss
queueing systems}

\author{
Vyacheslav M. Abramov\thanks{School of Mathematical Sciences,
Monash University, Building 28M, Clayton Campus, Clayton, VIC
3800, \textsc{Australia}.
\protect\url{Vyacheslav.Abramov@sci.monash.edu.au}}
\footnotemark[1] % when shares the address of 1st thanks
}

\date{18 April 2006}

\begin{document}

% Use default \verb|\maketitle|.
\maketitle

% Use the \verb|abstract| environment.
\begin{abstract}
%    The abstract is \emph{not} a table of contents.  Say what is
%    delivered, the essential qualities of the paper.  Use less than
%    50~words for each of the following questions: What was done?  Why
%    do it?  What were the results?  What do the results mean in theory
%    and/or practise?  What is the reader's benefit?  How can the
%    readers use this information for themselves?  The abstract is
%    probably all most readers read, it must be a complete description
%    in itself, though necessarily sketchy.

The present paper provides some new stochastic inequalities for
the characteristics of the $M/GI/1/n$ and $GI/M/1/n$ loss queueing
systems. These stochastic inequalities are based on substantially
deepen up- and down-crossings analysis, and they are stronger than
the known stochastic inequalities obtained earlier. Specifically,
for a class of $GI/M/1/n$ queueing system, two-side stochastic
inequalities are obtained.
\end{abstract}

% By default we include a table of contents in each paper.
\tableofcontents

% Use \verb|\section|, \verb|\subsection|, \verb|\subsubsection| and
% possibly \verb|\paragraph| to structure your document.  Make sure
% you \verb|\label| them for cross-referencing with \verb|\ref|\,.

\section{Introduction}\label{sec:intro}
The goal of the paper is to establish stronger stochastic
inequalities for the number of losses during a busy period than
those are obtained earlier in \cite{A3}. The number of losses
during a busy period is a significant characteristic for analysis
of loss probability and other performance measures of real
telecommunication systems, and detailed stochastic analysis of
losses in queueing systems seems to be very important.

For the purpose of detailed stochastic analysis of losses we
develop the up- and down-crossings approach initiated in a number
of earlier works of the author \cite{A1}-\cite{A4}. It is proved
in \cite{A3} that if the inter-arrival time distribution of
$GI/M/1/n$ queue belong to the class NBU (NWU), then the number of
losses during a busy period is stochastically not smaller
(respectively not greater) than the number of offspring in the
$n+1$st generation of the Galton-Watson branching process with
given offspring generating function (see below for the more
details). The Galton-Watson branching process is a well-known
process having relatively simple explicit expressions for its
characteristics. At the same time the explicit results for the
number of losses in the $M/GI/1/n$ and $GI/M/1/n$ queues are very
hard for applications.

In this paper we obtain two-side stochastic inequality  for the
number of losses during a busy period of the $GI/M/1/n$ queueing
system, where the left and right sizes are branching processes.

Note that other inequalities related to the number of losses
during a busy period in the different loss queueing systems were
obtained in \cite{PRX9}, \cite{R10}, \cite{W12} and others papers.

 The paper starts from elementary extension of the inequalities
obtained in \cite{A3} to some special class of $GI/GI/1/n$ queues,
which includes $M/GI/1/n$ queueing systems with NBU (NWU) service
time and $GI/M/1/n$ queueing systems with NBU (NWU) interarrival
time as particular cases. For our further convenience the
$GI/GI/1/n$ queueing system will be denoted $A/B/1/n$, where
$A(x)$ is the probability distribution function of an interarrival
time, and $B(x)$ is the probability distribution function of a
service time. Then, for the $M/GI/1/n$ and $GI/M/1/n$ queueing
system we often use the notation $M/B/1/n$ and $A/M/1/n$
respectively. For the definition of the classes of distributions
such as NBU, NWU and all other, that are used in the paper, see
\cite{S11}.

Throughout the paper the following notation is used. For $\Re
(s)\ge 0$ we denote  the Laplace-Stieltjes transforms of the
probability distributions $A(x)$ and $B(x)$ by $\widehat A(s)$ and
$\widehat B(s)$ respectively, and the reciprocals of the expected
inter-arrival and service times are denoted by $\lambda$ and $\mu$
respectively. The aforementioned Laplace-Stieltjes transforms are
in fact used for real values of argument, specifically only the
values $\widehat A(\mu)$ and $\widehat B(\lambda)$ are used
throughout the paper.

The number of losses during a busy period is denoted $L_n$.

For the $A/M/1/n$ queue we have the inequality
$L_n\ge_{st}Z_{n+1}$ in the case where an interarrival  time is
NBU, and the opposite inequality, $L_n\le_{st}Z_{n+1}$, in the
case where an interarrival time is NWU (see \cite{A3}). $Z_n$
denotes the number of offspring in the $n$th generation of the
Galton-Watson branching process with $Z_0=1$ and the offspring
generating function
\begin{equation*}
g_Z(z)=\frac{1-\widehat A(\mu)}{1-z\widehat A(\mu)}, \ |z|\le 1.
\end{equation*}
The method of \cite{A3}, adapted to the $M/B/1/n$ queue, provides
the following inequality:
\begin{equation}\label{1.1}
L_n\le_{st}Y_{n+1} \ \ \ \Big(L_n\ge_{st}Y_{n+1}\Big)
\end{equation}
in the case where the service time is NBU (NWU). $Y_n$ is the
number of offspring in the $n$th generation of the Galton-Watson
branching process with $Y_0=1$ and the offspring generating
function
\begin{equation*}
g_Y(z)=\frac{\widehat B(\lambda)}{1-z+z\widehat B(\lambda)}, \
|z|\le 1.
\end{equation*}
(See Section \ref{sect:2} for details of proof.)

A deepen analysis of these two queueing systems, given in Sections
\ref{sect:3} and \ref{sect:4}, enables us to obtain the following
stronger results than that permits us the method of \cite{A3}.

For the $M/B/1/n$ queue in the case where $B(x)$ belongs to the
class of NBU (NWU) distributions it is shown that
\begin{equation}\label{1.2}
L_n\le_{st}\sum_{i=1}^{Y_n}\tau_i \ \ \
\left(L_n\ge_{st}\sum_{i=1}^{Y_n}\tau_i\right),
\end{equation}
where $\tau_1, \tau_2,\ldots$ is a sequence of independent
identically distributed nonnegative integer random variables,
\[
\mathrm{P}[\tau_i=k]=\int_0^\infty\expon^{-\lambda
x}\frac{(\lambda x)^k}{k!}\deriv B(x).
\]

Representation \eqref{1.2} is preferable than \eqref{1.1}. For
example, it follows from \eqref{1.1} that
\begin{equation}\label{1.3}
\mathrm{E}[L_n]\ge\Big[\frac{1-\widehat B (\lambda)}{\widehat B
(\lambda)}\Big]^{n+1} \ \ \ \left(
\mathrm{E}[L_n]\le\Big[\frac{1-\widehat B(\lambda)}{\widehat
B(\lambda)}\Big]^{n+1}\right).
\end{equation}

In turn, by using the Wald's equation, from \eqref{1.2} we obtain
\begin{equation}\label{1.4}
\mathrm{E}[L_n]\ge\frac{\lambda}{\mu}\Big[\frac{1-\widehat B
(\lambda)}{\widehat B (\lambda)}\Big]^{n} \ \ \ \left(
\mathrm{E}[L_n]\le \frac{\lambda}{\mu}\Big[\frac{1-\widehat
B(\lambda)}{\widehat B(\lambda)}\Big]^{n}\right).
\end{equation}
Clearly that \eqref{1.4} is stronger than \eqref{1.3} since in the
case of the NBU (NWU) service time distribution we have:
\begin{equation*}
\frac{1-\widehat B(\lambda)}{\widehat
B(\lambda)}\geq\frac{\lambda}{\mu} \ \ \ \left(\frac{1-\widehat
B(\lambda)}{\widehat B(\lambda)}\leq\frac{\lambda}{\mu}\right).
\end{equation*}

For a subcritical $A/M/1/n$ queue ($\varrho=\lambda/\mu\le 1$), in
the case where an interarrival time distribution belongs to the
IHR (DHR) class of distributions we obtain $L_n\le_{st}X_{n+1}$
($L_n\ge_{st}X_{n+1}$) The process $\{X_n\}$ is a branching
process, but not classical (the precise definition of this process
is given in Section \ref{sect:4}). Thus, combining this result
with the result of \cite{A3} we conclude the following. If
$\varrho\le 1$, then in the case when an interarrival time
distribution belongs to the IHR (DHR) class of distributions we
have $Z_{n+1}\le_{st}L_n\le_{st}X_{n+1}$ \ \
($X_{n+1}\le_{st}L_n\le_{st}Z_{n+1}$).
\smallskip

The paper is organized as follows. It consists of 4 sections.
Section 2 introduces the reader to the up- and down-crossings
method of \cite{A3} and extends the results of \cite{A3} to the
special class of $A/B/1/n$ queues (described exactly in that
Section \ref{sect:2}). The results related to $M/B/1/n$ and
$A/M/1/n$ queues are then developed in Sections \ref{sect:3} and
\ref{sect:4} respectively. In turn, Section \ref{sect:4} is
divided into subsections, containing preliminary information on
the properties of the $A/M/1/n$ queues. The most significant
property is a monotonicity, which is considered in Section
\ref{sect:4.1}. Section \ref{sect:4.2} introduces and studies a
special type of branching process, which is then used for the main
result of Section \ref{sect:4} - Theorem \ref{th-4.3}.

\section{Stochastic inequalities for $GI/GI/1/n$ queues}\label{sect:2} In this section we
establish stochastic inequalities for a class of $A/B/1/n$ queues.
Specifically, assuming that the probability distributions $A(x)$
and $B(x)$ belong to the opposite classes of NBU and NWU, i.e
either $A(x)$ belongs to NBU and $B(x)$ belongs to NWU, or $A(x)$
belongs to NWU and $B(x)$ belongs to NBU, we have the following.

\medskip

    \begin{thm}
        \label{th-2.1}
 Under the assumption that $A(x)$ belongs to NBU
(NWU), and $B(x)$ belongs to NWU (NBU), and a busy period is
finite with probability 1, we have
\begin{equation}\label{2.1}
 L_n\ge_{st}X_{n+1} \ \ \ \Big(L_n\le_{st}X_{n+1}\Big).
\end{equation}
 $X_n$ in
\eqref{2.1} is the number of offspring in the $n$th generation of
the Galton-Watson branching process with $X_0=1$ and the offspring
generating function
\begin{equation}\label{2.2}
g(z)=\frac{1-r}{1-zr}, \ \ \ |z|\le 1,
\end{equation}
where $r=1-\int_0^\infty[1-A(x)]\deriv B(x)$.
    \end{thm}

\textit{Proof}. The proof is provided only in the case where
$A(x)$ belongs to the NBU class and $B(x)$ belongs to the NWU
class. The opposite case is analogous.

Let $f_n(j), 0\le j\le n+1$, denote the number of customers
arriving during a busy period who, upon their arrival, meet $j$
customers in the system. Under the assumption that a busy period
is finite we have $f_n(0)=1$ with probability 1. Let $t_{j,1}^n$,
$t_{j,2}^n$,\ldots, $t_{j,f_n(j)}^n$ be the instants of arrival of
these $f_n(j)$ customers, and let $s_{j,1}^n$, $s_{j,2}^n$,\ldots,
$s_{j, f_n(j)}^n$ be the instants of service completions
(departures) at which there remain only $j$ customers in the
system. Notice, that $t_{n+1,k}^n$ = $s_{n+1,k}^n$,~~ $1\leq k\leq
f_n(n+1)=L_n$.

For $0\le j\le n$ let us consider the following intervals:
\begin{equation}\label{2.3}
\Big[t_{j,1}^n, s_{j,1}^n\Big), ~\Big[t_{j,2}^n,
s_{j,2}^n\Big),..., \Big[t_{j,f_n(j)}^n, s_{j,f_n(j)}^n\Big).
\end{equation}
It is clear that the intervals
\begin{equation}\label{2.4}
\Big[t_{j+1,1}^n, s_{j+1,1}^n\Big), ~\Big[t_{j+1,2}^n,
s_{j+1,2}^n\Big),..., \Big[t_{j+1,f_n(j+1)}^n,
s_{j+1,f_n(j+1)}^n\Big)
\end{equation}
are contained in intervals \eqref{2.3}. Let us delete the
intervals in \eqref{2.4} from those in \eqref{2.3} and connect the
ends. That is, we connect every point $t_{j+1,k}^n$ with the
corresponding point $s_{j+1,k}^n,~1\le k\le f_n(j+1)$, if the set
of intervals \eqref{2.4} is not empty.

We will use the following notation. Take the interval $[t_{j,k}^n,
s_{j,k}^n)$. Within this interval there is a number of inserted
points, say $m$. If $m>0$ then these points are numbered as
$i=1,2,..., m$. Let $A_{j,k}^{(i)}(x)$ denote the probability
distribution of the residual time in point $i$ until the next
arrival, and let $B_{j,k}^{(i)}(x)$ denote  the probability
distribution of the residual service time in point $i$. Then
$A_{j,k}^{(0)}(x)$ is the probability distribution of the residual
time in the initial point $t_{j,k}^n$ of the interval $[t_{j,k}^n,
s_{j,k}^n)$ until the next arrival. Since $t_{j,k}^n$ is the
moment of arrival, then $A_{j,k}^{(0)}(x)=A(x)$ for all $j\ge 0$
and $k\ge 1$. $B_{j,k}^{(0)}(x)$ is the probability distribution
of the residual service time in the initial point $t_{j,k}^n$.

%A simple example of the up- and down-crossings on a busy period is
%given on Figure 1. The arc braces in the figure indicate the
%places of connection of one another points after deleting the
%intervals.

%\DeclareGraphicsRule{.ps}{ps}{.bb}{}%{.bb}{}
%\ \scalebox{0.8}{\includegraphics{figtex.ps}}

%In this figure the buffer length is assumed to be $n\geq 3$, the
%number of inserted points of the level 1 is equal to 2
%($f_n(1)=2$), and the number of inserted points of the level 2 is
%equal to 1 ($f_n(2)=1$).

\smallskip

Let us take the interval [$t_{j,k}^n, s_{j,k}^n$) and a customer
in service in time $t_{j,k}^n$. Let $\tau_{j,k}$ be the time
elapsed from the moment of the service begun for that customer
until time $t_{j,k}^n$. Then for residual service time
$\vartheta_{j,k}$ of the tagged customer we have
\begin{equation}\label{2.5}
\begin{aligned}
\mathrm{P}[\vartheta_{j,k}\le x]&=
\mathrm{P}[\chi\le\tau_{j,k}+x~|~\chi>\tau_{j,k}]\\
&= \int_0^\infty \mathrm{P}[\chi\le y+x~|~\chi>y]\deriv
\mathrm{P}[\tau_{j,k}\le y].
\end{aligned}
\end{equation}
According to the above convention, the probability of \eqref{2.5}
is denoted by $B_{j,k}^{(0)}(x)$. Let $\kappa_{j,k}$ denote the
number of inserted points within the interval [$t_{j,k}^n,
s_{j,k}^n$), so
\begin{equation*}
\sum_{k=1}^{f_n(j)}\kappa_{j,k}{\buildrel\Delta\over =}f_n(j+1).
\end{equation*}
Then,
\begin{equation*}
\mathrm{P}[\kappa_{j,k}=0]=\int_0^\infty[1-A_{j,k}^{(0)}(x)]\deriv
B_{j,k}^{(0)}(x),
\end{equation*}
and for $m\ge 1$
\begin{equation}\label{2.6}
\begin{aligned}
\mathrm{P}[\kappa_{j,k}=m]& =\prod_{i=0}^{m-1}\Big[1-\int_0^\infty
[1-A_{j,k}^{(i)}(x)]\deriv B_{j,k}^{(i)}(x)\Big]\times\\
&\times\int_0^\infty[1-A_{j,k}^{(m)}(x)]\deriv B_{j,k}^{(m)}(x).
\end{aligned}
\end{equation}
Relationship \eqref{2.6} looks cumbersome, but it has a simple
explanation. The term
\begin{equation*}
 \int_0^\infty[1-A_{j,k}^{(m)}(x)]\deriv B_{j,k}^{(m)}(x)
\end{equation*}
is the probability that during the residual service time
corresponding to the $m$th inserted point there is no arrival, or
in other words, the $m$th inserted point is last. Similarly, the
product term
\begin{equation*}
1-\int_0^\infty[1-A_{j,k}^{(i)}(x)]\deriv B_{j,k}^{(i)}(x)
\end{equation*}
is the probability that during the residual service time
corresponding to the $i$th inserted point there is at least one
arrival.

Taking into account that both $A(x)\leq A_{j,k}^{(i)}(x)$ and
$B_{j,k}^{(i)}(x)\leq B(x)$ for all $j,k$ and $i$, we have the
following. Let $\kappa_X$  be geometrically distributed random
variable, $\mathrm{P}[\kappa_X=m]=r^m(1-r)$, $m=0,1,...$, where
the parameters $r$  is determined in the formulation of the
theorem. Then, $\kappa_X\geq_{st}\kappa_{j,k}$, for all $j\ge 0$
and $k\ge 1$, and we have the following. Let $\kappa_X^{(j,k)}$ be
the sequences of independent identically distributed integer
random variables all having the same distribution as the random
variable $\kappa_X$. We have
\[
\sum_{k=1}^{f_n(j)}\kappa_{j,k}\leq_{st}\sum_{k=1}^{f_n(j)}\kappa_X^{(j,k)}.
\]
Taking into account that
\[
X_{j+1}=\sum_{k=1}^{X_j}\kappa_X^{(j,k)},
\]
owing to induction we have
\[
X_j\leq_{st}f_n(j),
\]
and therefore $f_n(n+1)=L_n\geq_{st}X_{n+1}$. The statement of the
theorem is proved. \ \ \ $\diamondsuit$

\smallskip

From Theorem \ref{th-2.1} we have the following special cases.

 \begin{cor}\label{cor-2.2}
 Under the assumption
that $A(x)=1-\expon^{-\lambda x}$, and $B(x)$ belongs to class NWU
(NBU), we have \eqref{2.1}. $X_n$ in \eqref{2.1} is the number of
offspring in the $n$th generation of the Galton-Watson branching
process with $X_0=1$ and the offspring generating function
\begin{equation}\label{2.7}
 g(z)=\frac{\widehat B(\lambda)}{1-z+z\widehat B(\lambda)}, \ \ \
|z|\le 1,
\end{equation}
 \end{cor}

\textit{Proof}. Putting $A(x)=1-\expon^{-\lambda x}$, we have
\begin{equation}\label{2.8}
r=1-\int_0^\infty \mbox{e}^{-\lambda x}\deriv B(x)=1-\widehat B
(\lambda),
\end{equation}
and the statement follows by substituting \eqref{2.8} for
\eqref{2.2}. \ \ \ $\diamondsuit$

\begin{cor}\label{cor-2.3} (Abramov [3].) Under the assumption
that $B(x)=1-\expon^{-\mu x}$, and $A(x)$ belongs to class NBU
(NWU), we have \eqref{2.1}. $X_n$ in \eqref{2.1} is the number of
offspring in the $n$th generation of the Galton-Watson branching
process with $X_0=1$ and the offspring generating function
\begin{equation*}
g(z)=\frac{1-\widehat A(\mu)}{1-z\widehat A(\mu)}, \ \ \ |z|\leq
1.
\end{equation*}
\end{cor}

\textit{Proof}. Putting $B(x)=1-\expon^{-\mu x}$, we have
\begin{equation}\label{2.9}
\begin{aligned}
r&=1-\int_0^\infty[1-A(x)]\mu\expon^{-\mu x}\deriv x\\
&=\int_0^\infty A(x)\mu\expon^{-\mu x}\deriv x\\
&=\int_0^\infty\expon^{-\mu x}\deriv A(x)\\ &=\widehat A(\mu).
\end{aligned}
\end{equation}
Substituting \eqref{2.9} for \eqref{2.2} we obtain the desired
representation. \ \ \ $\diamondsuit$

\section{Stronger inequalities for $M/GI/1/n$
queues}\label{sect:3} In this section we develop the result for
the $M/B/1/n$ queue given by Corollary \ref{cor-2.2}. The main
result of this section is the following.

\begin{thm}
Under the assumption that $A(x)=1-\expon^{-\lambda x}$, and $B(x)$
belongs to class NWU (NBU), we have
\begin{equation*}
L_n\ge_{st}\sum_{i=1}^{X_n}\tau_i \ \ \
\left(L_n\le_{st}\sum_{i=1}^{X_n}\tau_i\right),
\end{equation*}
where the branching process $\{X_n\}$ is the same as in Corollary
\ref{cor-2.2}, and $\tau_1, \tau_2,\ldots$ is a sequence of
independent identically distributed nonnegative integer random
variables,
\begin{equation*}
\mathrm{P}[\tau_i=k]=\int_0^\infty\expon^{-\lambda
x}\frac{(\lambda x)^k}{k!}\deriv B(x).
\end{equation*}

\end{thm}

\textit{Proof}.
 Considering first the $M/B/1/0$ loss queue without waiting
places it is not difficult to see that
\begin{equation*}
L_0~=_{st}~ \mbox{Number of Poisson arrivals per service time},
\end{equation*}
that is,
\begin{equation*}
\mathrm{P}[L_0=k]=\int_0^\infty\mbox{e}^{-\lambda x}\frac{(\lambda
x)^k}{k!}\deriv B(x).
\end{equation*}
Let us now consider the $M/B/1/n$ queue, where $f_n(n)$ is the
number of cases during a busy period when an arriving customer
meets $n$ customers in the system (recall that $L_n{\buildrel
\Delta\over =}f_n(n+1)$). Then, the number of losses $L_n$
coincides in distribution with
\begin{equation*}
\sum_{i=1}^{f_n(n)}\tau_i,
\end{equation*}
where the sequence $\tau_1, \tau_2,\ldots$ is a sequence of
independent identically distributed integer random variables,
coinciding in distribution with $L_0$.

It follows from the proof of Theorem \ref{th-2.1}, that if in the
$A/B/1/n$ queue $A(x)$ is NBU (NWU) and $B(x)$ is NWU (NBU), then
\begin{equation}\label{3.1}
f_n(n)\ge_{st} X_n \ \ \ \Big(f_n(n)\le_{st} X_n\Big)
\end{equation}
where the branching process $\{X_n\}$ is defined in Theorem
\ref{th-2.1}, i.e. $X_0=1$, and the offspring generating function
is determined by \ref{2.2}. Therefore, in the case of
$A(x)=1-\mbox{e}^{-\lambda x}$ we obtain \ref{3.1}, where now the
offspring generating function of the branching process is defined
by \eqref{2.7}. This enables us to conclude that under the
assumptions of the theorem
\begin{equation*}
L_n\ge_{st}\sum_{i=1}^{X_n}\tau_i \ \ \
\left(L_n\le_{st}\sum_{i=1}^{X_n}\tau_i\right),
\end{equation*}
and the statement is therefore proved.

Considering now the $A/B/1/n$ queueing system, let $T_n$, $\nu_n$
denote the length of a busy period and the number of served
customers during a busy period respectively,
 and let $\chi^{(1)}$, $\chi^{(2)}$,\ldots be a sequence
of independent identically distributed random variables all having
the probability distribution function $B(x)$. We have
\begin{equation*}
\nu_n=\sum_{j=0}^nf_n(j),
\end{equation*}
\begin{equation*}
T_n=\sum_{m=1}^{\nu_n}\chi^{(m)}.
\end{equation*}

Immediately from the above proof, under the assumption that $A(x)$
is NBU (NWU), and $B(x)$ is NWU (NBU), we have
\begin{equation}\label{3.2}
\nu_n\ge_{st}\sum_{i=0}^nX_i, \ \ \
\left(\nu_n\le_{st}\sum_{i=0}^nX_i\right),
\end{equation}
where the branching process $\{X_n\}$ is defined in Theorem
\ref{th-2.1}. If $A(x)=1-\mbox{e}^{-\lambda x}$, then \eqref{3.2}
holds true. The only difference that the offspring generating
function of the process $\{X_n\}$ is given by \eqref{2.7}.

Whereas the sequence of $\chi^{(1)}$, $\chi^{(2)}$,\ldots consists
of independent identically distributed random variables, the
random variable $\nu_n$ is independent of the future, that is the
event $\{\nu_n=i\}$ is independent of $\chi^{(i+1)}$,
$\chi^{(i+2)}$, \ldots (e.g. \cite{B6}). Therefore
$\mathrm{E}[T_n]$ is determined by the Wald's identity: $\mu
\mathrm{E}[T_n]=\mathrm{E}[\nu_n]$. Then under the above
assumptions that $A(x)$ is NBU (NWU) and $B(x)$ is NWU (NBU), we
have
\begin{equation}\label{3.3}
 \mathrm{E}[T_n]\geq\frac{1}{\mu}\mathrm{E}\Big[\sum_{i=0}^nX_i\Big] \ \ \
\left(\mathrm{E}[T_n]\leq\frac{1}{\mu}\mathrm{E}\Big[\sum_{i=0}^nX_i\Big]\right).
\end{equation}
 Taking into account that
$\mathrm{E}[X_n]=r^n/(1-r)^n$, under the above assumptions from
\eqref{3.3} we obtain
\begin{equation*}
\mathrm{E}[T_n]\geq\frac{1}{\mu}\sum_{i=0}^n\frac{r^i}{(1-r)^i} \
\ \
\left(\mathrm{E}[T_n]\leq\frac{1}{\mu}\sum_{i=0}^n\frac{r^i}{(1-r)^i}\right).
\end{equation*}
Clearly, that in the case where $A(x)=1-\expon^{-\lambda x}$ the
parameter $r$ is equal to $1-\widehat B(\lambda)$ (see the proof
of Corollary \ref{cor-2.2}).

\section{Further stochastic inequalities for the $GI/M/1/n$ loss
system}\label{sect:4} \noindent Being the special case of Theorem
\ref{th-2.1}, Corollary \ref{cor-2.3} provides the stochastic
inequalities for the $A/M/1/n$ under the assumption that $A(x)$
belongs to the class NBU (NWU). Assuming now that $A(x)$ belongs
to the class IHR (DHR), we provide a deepen analysis in order to
obtain stronger stochastic inequalities.

\subsection{Monotonicity}\label{sect:4.1}
For the sake of simplicity the $A/M/1/n$ queueing system is
denoted $\mathcal{Q}_n$. Recall that parameter $n$ excludes the
position of a customer in service. For $n$ and $k$ different,
$\mathcal{Q}_n$, $\mathcal{Q}_k$ are two queueing systems with the
same probability distribution functions of interarrival and
service time but different number of waiting places. For example,
$\mathcal{Q}_0$ denotes the $A/M/1/0$ queueing system without
waiting places, a busy period of which contains only a single
service time.

Consider a busy period of the queueing system $\mathcal{Q}_n$. Let
us consider the interval $[t_{0,1}^n, s_{0,1}^n)$ after the
procedure of deleting from it all the intervals $[t_{1,l}^n,
s_{1,l}^n)$, ~$l=1,2,...,f_n(1)$, and connecting the ends as it is
described in the proof of Theorem \ref{th-2.1}. Then, let $i_1^n$,
$i_2^n,...$, $i_{f_n(1)}^n$ denote the inserted points within the
interval $[t_{0,1}^n, s_{0,1}^n)$, and let $d_l^n$ denote the
distance between the two adjacent points $i_{l}^n$ and $i_{l+1}^n$
~($l=1,2,...,f_n(1)-1$), that is, $d_l^n=i_{l+1}^n-i_{l}^n$. If
$f_n(1)=0$, i.e. there is no inserted points, then the distance
between inserted points is not defined. If $f_n(1)=1$, then by the
value $d_1^n$ we mean the distance between the point $i_1^n$ and
the next arrival of a customer at the system after the instant
$s_{0,1}^n$.

\medskip

\begin{lem}\label{lem-4.1} Let $\mathcal{Q}_k$ and $\mathcal{Q}_n$
be two queueing systems, and let $A(x)$ belong to the IHR (DHR)
class of distributions. If $k\le n$ then
\begin{equation*}
d_l^n\le_{st}d_l^k \ \ \ \Big(d_l^n\ge_{st}d_l^k\Big).
\end{equation*}
\end{lem}

\textit{Proof}. Let us consider the queueing system
$\mathcal{Q}_n$, and the interval $[t_{0,1}^n$, $s_{0,1}^n)$ after
the procedure of deleting from it all intervals $[t_{1,l}^n,
s_{1,l}^n)$, ~$l=1,2,...$, $f_n(1)$, and connecting the ends. For
convenience, we denote the sequence of independent identically and
exponentially distributed random variables with parameter $\mu$ by
$\chi^{(1)}$, $\chi^{(2)},...$, and a random variable $\tau$,
having the probability distribution $A(x)$, is independent of this
sequence $\chi^{(1)}$, $\chi^{(2)},\ldots$.

The probability, that during the interval $[t_{0,1}^n, s_{0,1}^n)$
there is no arrival, is
\begin{equation*}
1-\int_0^\infty\mu\mbox{e}^{-\mu x}A(x)\deriv x=1-\widehat A(\mu).
\end{equation*}
Obviously, that this probability is independent of parameter $n$.
Let us assume that there is the inserted point $i_1^n$ and,
therefore, the instant of arrival $t_{1,1}^n$.

Let $q_n$ denote the stationary number of customers in the
queueing system $\mathcal{Q}_n$ immediately after arrival of a
customer at the system \emph{during a busy period}, i.e. not into
the empty system. (An arriving customer, who finds all waiting
places busy, leaves the system without incrementing and
decrementing the number of customers in the queue.) Let
$\widetilde q_n$ = $q_n-1$, and let
\begin{equation*}
v=\inf\Big\{ m: \sum_{j=1}^{\widetilde
q_n}\chi^{(j+m-1)}\le\tau\Big\}.
\end{equation*}
(The empty sum is assumed to be 0. The case of empty sum arises
only by considering of the queueing system $\mathcal{Q}_0$.) Then
\begin{equation}\label{4.1}
d_l^n~{\buildrel 'd'\over =}~\tau-\sum_{j=1}^{\widetilde
q_n}\chi^{(v+j-1)}.
\end{equation}
For example, in the case of the queueing system $\mathcal{Q}_0$,
we have
\begin{equation*}
\mathrm{P}[d_{l}^0\le x]=\mathrm{P}[\tau\le x]=A(x),
\end{equation*}
and in the case of the queueing system $\mathcal{Q}_1$ we have
\begin{equation}\label{4.2}
\begin{aligned}
\mathrm{P}[d_{l}^1\le x]&=\mathrm{P}[ \tau - \chi^{(1)}\le x|
\tau>\chi^{(1)}]\\
&= \int_0^\infty \mathrm{P}[\tau\le x+y|\tau>y]\mu\expon^{-\mu
y}\deriv y\\
&=\int_0^\infty\frac{A(x+y)-A(y)}{1-A(y)}\mu\expon^{-\mu y}\deriv
y.
\end{aligned}
\end{equation}
By analysis of sample paths it is clear that for these two
queueing systems $\mathcal{Q}_n$ and $\mathcal{Q}_{n+1}$
\begin{equation}\label{4.3}
\widetilde q_n\leq_{st}\widetilde q_{n+1}.
\end{equation}
Since $A(x)$ belongs to the IHR (DHR) class of distributions, then
\eqref{4.3} together with \eqref{4.1} yield
$d_l^{n+1}\leq_{st}d_l^n$ \ ($d_l^{n+1}\geq_{st}d_l^n$). The
statement of lemma follows. \ \ \ $\diamondsuit$

\begin{rem} Lemma \ref{lem-4.1} establishes a property of external
monotonicity. However, from Lemma \ref{lem-4.1} we obtain the
property of internal monotonicity as well. Indeed, in the case of
the $GI/M/1/n$ queueing system, because of the property of the
lack of memory of the exponential distribution of a service time,
any interval of \eqref{2.3} is distributed as a busy period of the
queueing system $\mathcal{Q}_{n-j}$, \ $0\le j\le n$. Therefore
the distance between two inserted points of each interval
\eqref{2.3} coincides in distribution with $d_1^{n-j}$.
\end{rem}

\subsection{A branching process}\label{sect:4.2} \noindent Let us consider the
$A/M/1$ queueing system (with infinite number of waiting places),
denoting it by $\mathcal{Q}$ and remaining for this system all the
above notation  given earlier for the queueing system
$\mathcal{Q}_n$. Assume additionally that the load
$\rho=\lambda/\mu\le 1$.

Analogously to the case of the queueing system $\mathcal{Q}_n$,
for the queueing system $\mathcal{Q}$ let  $f(j)$,~ $j\ge 0$,
denote the number of customers, arriving during a busy period,
who, upon their arrival meet $j$ customers in the system
($f(0)=1$). Let $t_{j,1}$, $t_{j,2},...$, $t_{j,f(j)}$ be the
instants of these arrivals, and let $s_{j,1}$, $s_{j,2},...$,
$s_{j,f(j)}$ be the instants of corresponding service completions
defined analogously to the case of the queueing system
$\mathcal{Q}_n$. Let $\mathcal{F}_j$ = $\sigma\{f(0),
f(1),...,f(j)\}$.

It is claimed in \cite{A4}, that the stochastic sequence $\{ f(j),
\mathcal{F}_j\}$ is a Galton-Watson branching process, and $
\mathrm{E}[f(1)]=\varphi$, where $\varphi$ is the least in
absolute value root of the functional equation $z=\widehat
A(\mu-\mu z)$.

According to the standard definition of the Galton-Watson
branching process, the number of offspring generated by all
particles are mutually independent random variables (e.g. Harris
\cite{H8}). The Galton-Watson branching process $\{ f(j),
\mathcal{F}_j\}$, considered in \cite{A4} for the case of $GI/M/1$
queues, is not traditional. The number of offspring generated by
particles of different generations are not independent random
variables. More precisely, the number of offspring of the $n$th
generation is an independent of the future random variable with
respect to the numbers of offspring generated by particles of the
$n$th generation.

Notice, that connection between standard branching process and
optimal stopping times has been discussed by Assaf, Goldstein and
Samuel-Sahn \cite{AGS5}.

For a more detailed explanation the structure of the
abovementioned dependence, related to the above case of the
$A/M/1$ queueing system, let us consider the interval [$t_{0,1},
s_{0,1}$), and assume that there is a point $t_{1,1}$. Let
$d_1=t_{1,2}-s_{1,1}$ denote the distance between the begin of the
second interval and the end of the first one (provided that the
second interval does exist). If there is only a single interval
then $d_1$ also has sense as it is explained in Section
\ref{sect:4.1}.

If during the time interval [$t_{1,1}, s_{1,1}$) there is no new
arrival (denote this event by $E_0$), then
\begin{equation}\label{4.4}
\begin{aligned}
\mathrm{P}[d_1\le x|E_0]&=\mathrm{P}[ \tau - \chi_1\le x|
\tau>\chi_1]\\ &= \int_0^\infty \mathrm{P}[\tau\le
x+y|\tau>y]\mu\expon^{-\mu y}\deriv y\\
&=\int_0^\infty\frac{A(x+y)-A(y)}{1-A(y)}\mu\expon^{-\mu y}\deriv
y.
\end{aligned}
\end{equation}
Recall that $\mathrm{P}[\tau\le x]=A(x)$, and
$\mathrm{P}[\chi_1\le x]=1-\expon^{-\mu x}$. Thus \eqref{4.4}
coincides with \eqref{4.2}, and $ \mathrm{P}[d_1\le x|E_0]$ =
$\mathrm{P}[d_1^1\leq x]$. For example, if $ \mathrm{P}[\tau=1]=1$
and $\mu\ge 1$, then from \eqref{4.4} we obtain
\begin{equation*}
\mathrm{P}[d_1\leq x|E_0]=\min\Big\{1, \frac{\mbox{e}^{\mu
x-\mu}-\expon^{-\mu}}{1-\expon^{-\mu}}\Big\}, \ \ x\ge 0.
\end{equation*}
If during the time interval [$t_{1,1}, s_{1,1}$) there is at least
one arrival (denote this event by $E_1$), then we have the
following. Let $\{q(i)\}_{i\ge 1}$ be a stationary sequence of the
numbers of customers in the system immediately {\it before}
arrival of a customer during a busy period (i.e. not into the
empty system). Let us consider the sequence $\{q(i){\rm
1}_{[q(i)\ge 2]}\}_{i\ge 1}$. Taking only the positive elements of
this sequence one can construct a {\it new stationary sequence}
$\{\widetilde q(i)\}_{i\ge 1}$ all elements of which are not
smaller than 2. Then,
\begin{equation*}
v=\inf\Big\{ m: \sum_{j=1}^{\widetilde
q(1)}\chi^{(j+m-1)}\leq\tau\Big\},
\end{equation*}
and
\begin{equation}\label{4.5}
\mathrm{P}[d_1\leq x|E_1]= \mathrm{P}\Big[\tau-\sum_{j=1}^{
\widetilde q(1)}\chi^{(j+v-1)}\leq x\Big].
\end{equation}
Comparing \eqref{4.4} and \eqref{4.5} it is not difficult to
conclude that if $A(x)$ belongs to the IHR (DHR) class of
distributions, then
\begin{equation*}
\begin{aligned}
&\mathrm{P}[d_1\le x|E_0]\leq \mathrm{P}[d_1\leq x|E_1]\\ \Big(
&\mathrm{P}[d_1\le x|E_0]\geq \mathrm{P}[d_1\leq x|E_1]\Big)
\end{aligned}
\end{equation*}
For example, if $\mathrm{P}[\tau=1]=1$, and $\mu\ge 1$, then we
have the strong inequality:
\begin{equation*}
\begin{aligned}
\mathrm{P}[d_1\leq x|E_0] &= \min\Big\{1, \frac{\expon^{\mu
x-\mu}-\expon^{-\mu}}{1-\expon^{-\mu}}\Big\}\\ &<
\mathrm{P}[d_1\leq x|E_1] \ \ \ (x\ge 0).
\end{aligned}
\end{equation*}
Thus, the random variable $f(1)$ depends on the events $E_0$ and
$E_1$. In other words $f(1)$ can have different distributions if a
particle of the first generation has or does not have an
offspring. Let us call such Galton-Watson branching process by
$GI/M/1$ {\it type Galton-Watson branching process}.

Notice, that the known property of a Galton-Watson branching
process that $\mathrm{E}[f(j)]=\varphi^j$ (e.g. Doob \cite{D7},
Harris \cite{H8}), is also remain in force for the $GI/M/1$ type
Galton-Watson branching process.

Indeed, according to the total expectation formula, for
$\mathrm{E}[f(1)]$ we obtain:

%--------------------- Original version ---------------------------
%
%\begin{equation}\label{r4.13}
% E[f(1)]=\sum_{n=0}^\infty E[f(n)]\int_0^\infty\mbox{e}^{-\mu
%x}\frac{(-\mu x)^n}{n!}\mbox{d}A(x).
%\end{equation}
%By the same arguments for all $j\geq 1$ we have:
%\[
%E[f(j+1)]=\sum_{n=0}^\infty E[f(n+j)]\int_0^\infty\mbox{e}^{-\mu
%x}\frac{(-\mu x)^n}{n!}\mbox{d}A(x).
%\]
%Therefore $Ef(n)=z^n$, and from (\ref{r4.13}) we have:
%\[
% E[f(1)]=z=\sum_{n=0}^\infty z^n\int_0^\infty\mbox{e}^{-\mu
%x}\frac{(-\mu x)^n}{n!}\mbox{d}A(x)=\widehat A(\mu-\mu z).
%\]
%Since $z<1$, then $z=\varphi$, and $E[f(n)]=\varphi^n$.

%----------------- Correction: 16 April 2005 ------------------

\begin{equation}\label{4.6}
 \mathrm{E}[f(1)]=\sum_{n=0}^\infty \mathrm{E}[f(n)]\int_0^\infty\expon^{-\mu
x}\frac{(\mu x)^n}{n!}\deriv A(x)
\end{equation}
By
the same arguments for all $j\geq 1$ we have:
\begin{equation*}
\mathrm{E}[f(j+1)]=\sum_{n=0}^\infty
\mathrm{E}[f(n+j)]\int_0^\infty\expon^{-\mu x}\frac{(\mu
x)^n}{n!}\deriv A(x).
\end{equation*}
Therefore $\mathrm{E}f(n)=z^n$, and from \eqref{4.6} we have:
\begin{equation*}
 \mathrm{E}[f(1)]=z=\sum_{n=0}^\infty z^n\int_0^\infty\expon^{-\mu
x}\frac{(\mu x)^n}{n!}\deriv A(x)=\widehat A(\mu-\mu z).
\end{equation*}
Since $z<1$, then $z=\varphi$, and $\mathrm{E}[f(n)]=\varphi^n$.

\subsection{The number of losses
during a busy period}\label{sect:4.3} Returning to the queueing
system $\mathcal{Q}_n$ once again, assume additionally that the
load $\rho\le 1$. All queueing systems $\mathcal{Q}_n$ with
different $n$ and the queueing system $\mathcal{Q}$ are assumed to
be given on the same probability space, and the probability
distribution function $A(x)$ belongs to the IHR (DHR) class of
distributions. According to Lemma \ref{lem-4.1} we have
\begin{equation}\label{4.7}
 d_l\le_{st}d_l^n \ \ (d_l\ge_{st}d_l^n),
\end{equation}
where $d_l$ is the distance between the $l$th and $l+1$st inserted
points of the queueing system $\mathcal{Q}$, as it is precisely
defined in Section \ref{sect:4.2}. Stochastic inequality
\eqref{4.7} is the limiting case, as $k\to\infty$, of a series of
inequalities for the distances $d_l^k\leq_{st}d_l^n$
($d_l^k\geq_{st}d_l^n$), given for all $k>n$.

Let us now consider the interval $[t_{0,1}^n, s_{0,1}^n)$ after
deleting all the intervals $[t_{1,j}^n, s_{1,j}^n)$ and connecting
the ends, as it is explained above. Then the remaining interval,
because of the property of the lack of memory, is exponentially
distributed with parameter $\mu$, and it coincides in distribution
with the  interval $[t_{0,1}, s_{0,1})$, associated with the
queueing system $\mathcal{Q}$, remaining after deleting of all the
intervals $[t_{1,j}, s_{1,j})$ and connecting the ends. Under the
assumption that both queueing processes of $\mathcal{Q}_n$ and
$\mathcal{Q}$ are defined on the same probability space, one may
consider only one of these intervals, comparing then the sample
path of relevant processes. Then for the number of losses $L_n$
during a busy period of the queueing system $\mathcal{Q}_n$ we
have the following.

\begin{thm}\label{th-4.3}
 If $A(x)$ belongs to the IHR (DHR)
class of distributions, and the load $\rho\le 1$, then
\begin{equation*}
 L_n\le_{st}Y_{n+1} \ \ \ \Big(L_n\ge_{st}Y_{n+1}\Big),
\end{equation*}
where $Y_n$ denotes the number of offspring in the $n$th
generation of the $GI/M/1$ type Galton-Watson branching process
generated by the queueing system $\mathcal{Q}$.
\end{thm}

Notice, that under the assumptions of Theorem \ref{th-4.3} we have
the inequality
\begin{equation}\label{4.8}
\mathrm{E}[L_n]\leq\varphi^{n+1}~ ~
~\Big(\mathrm{E}[L_n]\ge\varphi^{n+1}\Big).
\end{equation}
On the other hand, taking into account that the class
 IHR (DHR) is contained in the class  NBU (NWU), from Corollary \ref{cor-2.3}
 we obtain the
inequality:
\begin{equation}\label{4.9}
\begin{aligned}
&\Big[\frac{\widehat A(\mu)}{1-\widehat A(\mu)}\Big]^{n+1}\leq
\mathrm{E}[L_n]\\
&\left(\Big[\frac{\widehat A(\mu)}{1-\widehat A
(\mu)}\Big]^{n+1}\geq \mathrm{E}[L_n]\right).
\end{aligned}
\end{equation}
Joining \eqref{4.8} and \eqref{4.9}, under the assumptions of
Theorem \ref{th-4.3} we obtain the two-side inequalities
\begin{equation}\label{4.10}
\begin{aligned}
&\Big[\frac{\widehat A(\mu)}{1-\widehat A(\mu)}\Big]^{n+1}\leq
\mathrm{E}[L_n]\le\varphi^{n+1}\\
&\left(\Big[\frac{\widehat A(\mu)}{1-\widehat
A(\mu)}\Big]^{n+1}\ge E[L_n]\ge\varphi^{n+1}\right).
\end{aligned}
\end{equation}
For example, in the case of the $M/M/1/n$ queueing system, when
$A(x)=1-\expon^{-\lambda x}$, from \eqref{4.10} we obtain $
\mathrm{E}[L_n]=\rho^{n+1}$.

\smallskip
It is interesting to note the following property. It follows from
\eqref{4.8} that if $A(x)$ belongs to the IHR (DHR) class of
distributions and $\rho\le 1$ ($\rho\ge 1$), then $
\mathrm{E}[L_n]\le 1$ ($ \mathrm{E}[L_n]\ge 1$) for all $n\ge 0$.
This is the special case of the more general result of Wolff
\cite{W12} for losses in $GI/GI/1/n$ queues under the assumption
that interarrival time distribution belongs to the class NBUE
(NWUE).

\smallskip
Let us provide inequalities for a busy period $T_n$ and the number
of customers served during a busy period of the $A/M/1/n$ queue.
Under the assumption that $A(x)$ is IHR (DHR) and $\rho<1$, we
have
\begin{equation}\label{4.11}
 \nu_n\leq_{st}\sum_{j=0}^nY_j \ \ \
\Big(\nu_n\geq_{st}\sum_{j=0}^nY_j \Big),
\end{equation}
where the branching process $\{Y_j\}$ is as in Theorem
\ref{th-4.3}.

From \eqref{4.11}, assuming that $A(x)$ is IHR (DHR) and $\rho<1$,
we obtain
\begin{equation}\label{4.12}
\mathrm{E}[\nu_n]\le\sum_{i=0}^n\varphi^{i} \ \ \ \left(
\mathrm{E}[\nu_n]\ge\sum_{i=0}^n\varphi^{i}\right).
\end{equation}
On the other hand, taking into account that class
 IHR (DHR) is contained in class  NBU (NWU), from Corollary \ref{cor-2.3} we obtain the
 following inequality:
\begin{equation}\label{4.13}
\begin{aligned}
&\sum_{i=0}^n\Big[\frac{\widehat A(\mu)}{1-\widehat
A(\mu)}\Big]^{i}\leq \mathrm{E}[\nu_n]\\
&\left(\sum_{i=0}^n\Big[\frac{\widehat A(\mu)}{1-\widehat
A(\mu)}\Big]^{i}\ge \mathrm{E}[\nu_n]\right)
\end{aligned}
\end{equation}
Combining \eqref{4.12} and \eqref{4.13}, under the above
assumptions we obtain the two-side inequalities:
\begin{equation*}
\begin{aligned}
&\sum_{i=0}^n\Big[\frac{\widehat A(\mu)}{1-\widehat
A(\mu)}\Big]^{i}\leq \mathrm{E}[\nu_n]\leq\sum_{i=0}^n\varphi^{i}\\
&\left(\sum_{i=0}^n\Big[\frac{\widehat A(\mu)}{1-\widehat
A(\mu)}\Big]^{i}\geq
\mathrm{E}[\nu_n]\geq\sum_{i=0}^n\varphi^{i}\right).
\end{aligned}
\end{equation*}
Finally, by Wald's identity we have
\begin{equation*}
\begin{aligned}
&\frac{1}{\mu}\sum_{i=0}^n\Big[\frac{\widehat A(\mu)}{1-\widehat
A(\mu)}\Big]^{i}\le
\mathrm{E}[T_n]\le\frac{1}{\mu}\sum_{i=0}^n\varphi^{i}\\
&\left(\frac{1}{\mu}\sum_{i=0}^n\Big[\frac{\widehat
A(\mu)}{1-\widehat A(\mu)}\Big]^{i}\geq
\mathrm{E}[T_n]\geq\frac{1}{\mu}\sum_{i=0}^n\varphi^{i}\right).
\end{aligned}
\end{equation*}

\section*{Acknowledgement}
The author thanks Taleb Samira (USTHB, Algeria) for useful
comment.


\begin{thebibliography}{99}
\bibitem{A1}
Abramov, V.M., {\em Investigation of a Queueing System with
Service Depending on Queue Length}. Donish, Dushanbe,
Tadzhikistan, 1991. (In Russian.)

\bibitem{A2}
Abramov, V.M., On the asymptotic distribution of the maximum
number of infectives in epidemic models with immigration, {\em
Journal of Applied Probability} {\bf 31} (1994), 606-613.

\bibitem{A3}
Abramov, V.M., Inequalities for the $GI/M/1/n$ loss system, {\em
Journal of Applied Probability} {\bf 38} (2001), 232-234.


\bibitem{A4}
Abramov, V.M., Some results for large closed queueing networks
with and without bottleneck: Up- and down-crossings approach, {\em
Queueing Systems} {\bf 38} (2001), 149-184.

%\bibitem{A5}
%Abramov, V.M., Asymptotic analysis of the $GI/M/1/n$ loss system
%as $n$ increases to infinity, \emph{Annals of Operations Research}
%{\bf 112} (2002) 35-41.
%
%
%
%
%\bibitem{A6}
%Abramov, V.M., Asymptotic behavior of the number of lost messages,
%\emph{SIAM Journal on Applied Mathematics} {\bf 64} (2004),
%746-761.


\bibitem{AGS5}
Assaf, D., Goldstein, L. and Samuel-Sahn, E., An unexpected
connection between branching processes and optimal stopping, {\em
Journal of Applied Probability} {\bf 37} (2000), 613-626.



\bibitem{B6}
Borovkov, A.A., {\em Theory of Probability.} Nauka, Moscow, 1986.
(In Russian.)


\bibitem{D7}
Doob, J.L., {\em Stochastic Processes.} John Wiley, New York,
1953.


\bibitem{H8}
Harris, T.E., {\em The Theory of Branching Processes}.
Springer-Verlag, Berlin, 1963.


\bibitem{PRX9}
Pek\"oz, E.A., Righter, R. and Xia, C.H., Characterizing losses in
finite buffer systems, {\em Journal of Applied Probability} {\bf
40} (2003), 242-249.



\bibitem{R10}
Righter, R., A note on losses in the $M/GI/1/n$ queue, {\em
Journal of Applied Probability} {\bf 36} (1999), 1240-1243.


\bibitem{S11}
Stoyan, D. {\em Comparison Methods for Queues and Other Stochastic
Models}. John Wiley, Chichester, 1983.

%\bibitem{S}
%Subhankulov, M.A., \emph{Tauberian Theorems with Remainder}.
%Nauka, Moscow, 1976. (In Russian.)

\bibitem{W12}
Wolff, R.W., Losses per cycle in a single-server queue, {\em
Journal of Applied Probability} {\bf 39} (2002), 905-909.
    \end{thebibliography}
\end{document}